\documentclass[12pt]{article}
\usepackage{amssymb}
\usepackage{amsmath}
\newcommand{\la}{\lambda}
\newcommand{\lap}{\mbox{$\bigtriangleup$}}

\newcommand{\ra}{{\mbox{$\rightarrow$}}}
\newcommand{\be}{\begin{equation}}
\newcommand{\ee}{\end{equation}}

\newtheorem{mthm}{Theorem}

\newtheorem{pro}{Proposition}[section]

\begin{document}
\title{A direct blowing-up and rescaling argument on the fractional Laplacian equation}

\author{Wenxiong Chen \thanks{Partially supported by the Simons Foundation Collaboration Grant for Mathematicians 245486.} \quad  Congming Li \thanks{ Corresponding author, partially supported by NSF DMS-1405175 and NSFC-11271166.} \quad Yan Li  }

\date{\today}
\maketitle
\begin{abstract} In this paper, we develop a direct {\em blowing-up and rescaling} argument for a nonlinear equation involving the fractional Laplacian operator. Instead of using the conventional extension method introduced by Caffarelli and Silvestre, we work directly on the nonlocal operator. Using the integral defining the nonlocal elliptic operator, by an elementary approach, we carry on a {\em blowing-up and rescaling} argument directly on nonlocal equations and thus obtain a priori estimates on the positive solutions for a semi-linear equation involving the fractional Laplacian.

We believe that the ideas introduced here can be applied to problems involving more general nonlocal operators.
\end{abstract}
\bigskip

{\bf Key words} The fractional Laplacian, blowing-up, rescaling, a priori estimates.
\bigskip

\section{Introduction}

The fractional Laplacian in $R^n$ is a nonlocal pseudo-differential operator, assuming the form
\begin{eqnarray}
(-\Delta)^{\alpha/2} u(x)& = &C_{n,\alpha} \, \lim_{\epsilon \ra 0} \int_{\mathbb{R}^n\setminus B_{\epsilon}(x)} \frac{u(x)-u(z)}{|x-z|^{n+\alpha}} dz\\
& = &C_{n,\alpha}\,PV\int_{\mathbb{R}^n} \frac{u(x)-u(z)}{|x-z|^{n+\alpha}} dz,
\label{op}
\end{eqnarray}
 where $\alpha$ is any real number between $0$ and $2$.
This operator is well defined in $\cal{S}$, the Schwartz space of rapidly decreasing $C^{\infty}$
functions in $\mathbb{R}^n$. In this space, it can also be equivalently defined in terms of the Fourier transform
$$ \widehat{(-\Delta)^{\alpha/2} u} (\xi) = |\xi|^{\alpha} \hat{u} (\xi), $$
where $\hat{u}$ is the Fourier transform of $u$. One can extend this operator to a wider space of functions.

Let
$$L_{\alpha}=\{u: \mathbb{R}^n\rightarrow \mathbb{R} \mid \int_{\mathbb{R}^n}\frac{|u(x)|}{1+|x|^{n+\alpha}} \, d x <\infty\}.$$
Then it is easy to verify that for $u \in L_{\alpha}  \cap C_{loc}^{1,1}$, the integral on the right-hand side of (\ref{op}) is well defined. Throughout this paper, we consider the fractional Laplacian in this setting.

The non-locality of the fractional Laplacian makes it difficult to investigate. To circumvent this difficulty,
Caffarelli and Silvestre \cite{CS} introduced the {\em extension method} that reduced a nonlocal problem into a local one in higher dimensions. For a function $u:\mathbb{R}^n \ra \mathbb{R}$, consider the extension $U:\mathbb{R}^n\times[0, \infty) \ra \mathbb{R}$ that satisfies
$$\left\{\begin{array}{ll}
div(y^{1-\alpha} \nabla U)=0, & (x,y) \in \mathbb{R}^n\times[0, \infty),\\
U(x, 0) = u(x).
\end{array}
\right.
$$
Then
$$(-\lap)^{\alpha/2}u (x) = - C_{n,\alpha} \displaystyle\lim_{y \ra 0^+}
y^{1-\alpha} \frac{\partial U}{\partial y},  \;\; x \in \mathbb{R}^n.$$

This {\em extension method } has been applied successfully to problems involving the fractional Laplacian, and a series of fruitful results have been obtained (see \cite{BCPS} \cite{BCPS1} \cite{CZ}    and the references therein).

The main purpose of this paper is to introduce a direct approach to study these nonlocal operators. We will carry on the blowing up and rescaling arguments on the nonlocal problem directly to obtain the a priori estimates on the solutions and to obtain the following.

\begin{mthm}
For $0<\alpha<2$ and a bounded domain $\Omega \subset \mathbb{R}^n$, consider
\begin{equation}
\left\{\begin{array}{ll}
(-\lap)^{\alpha/2}u(x)=u^p(x),  & x \in \Omega,\\
u(x)\equiv0, &x \not\in \Omega.
\end{array}\label{pe1}
\right.
\end{equation}
Assume that $u\in L_\alpha(\mathbb{R}^n)\cap C^{1,1}_{loc}(\Omega)$ and $1<p<\frac{n+\alpha}{n-\alpha}$, then there exists some
constant $C$, independent of $u$, such that
\be
\parallel u\parallel_{L^\infty(\Omega)}\leq C.\label{pe2}
\ee
\label{mthm1}
\end{mthm}

It is well-known that a priori estimates play important roles in establishing the existence of solutions.
Once there is such an a priori estimate, then one can use various methods, such as continuation or topological degrees, to derive the existence of solutions.

\section{Proof of Theorem \ref{mthm1}}

\textbf{Proof.} Suppose (\ref{pe2}) does not hold, then there exists a sequence
 of solutions $\{u_k\}$ to (\ref{pe1}) and a sequence of points $\{x^k\} \subset \Omega$ such that
\be
u_k(x^k)=\underset{\Omega}{\max}\,u_k:=m_k \ra \infty.
\ee

Let
\be
\lambda_k=m_k^{\frac{1-p}{\alpha}} \mbox{ and } v_k(x)=\frac{1}{m_k}u_k(\lambda_k x+x^k),\label{pe2.1}
\ee
then we have
\be
(-\lap)^{\alpha/2}v_k(x)=v_k^p(x), \quad x
 \in \Omega_k:=\{x \in \mathbb{R}^n \mid x=\frac{y-x^k}{\lambda_k},\,
y \in \Omega\}. \label{pe3}
\ee
Let $d_k=dist(x^k, \partial\Omega)$.
We will carry out the proof using the contradiction argument while
exhausting all three possibilities.
\vspace{12pt}

\emph{Case i.} $\underset{k \ra \infty}{\lim}\frac{d_k}{\lambda_k} =\infty$.
\vspace{12pt}

It's not difficult to see that
$$\Omega_k \ra \mathbb{R}^n \mbox{ as } k \ra \infty.$$
 We'll prove that there exists a function $v$ such that as
 $k \ra \infty$,
\be
v_k(x) \ra v(x) \mbox{ and }
 (-\lap)^{\alpha/2}v_k(x)\ra(-\lap)^{\alpha/2}v(x),\label{pe4}
\ee
thus
\be
(-\lap)^{\alpha/2}v(x)=v^p(x), \quad x \in  \mathbb{R}^n. \label{pe5}
\ee

It has been proved that (\ref{pe5})
has no positive solution (see ). Meanwhile, (\ref{pe2.1}) indicates that
$$v(0)=\lim_{k \ra \infty} v_k(0)=1.$$
This is a contradiction. Hence $u$ must be uniformly bounded in $\Omega$.

To obtain (\ref{pe4}), we need the following propositions to boost the regularity
of $v_k$.

\begin{pro}
\cite{Si}
Let $u\in C^{k, \sigma}$ and suppose that $k+\sigma-\alpha$ is not an integer.
 Then
$(-\lap)^{\alpha/2} u \in C^{l,\gamma}$ where $l$ is the integer part of
$k+\sigma-\alpha$ and $\gamma = k+\sigma-\alpha-l$.
\label{prop2.7}
\end{pro}

\begin{pro}
\cite{Si}
Let $w=(-\lap)^{\alpha/2}u$. Assume that $w \in L^\infty(\mathbb{R}^n)$ and
$u \in L^\infty(\mathbb{R}^n)$ for $\alpha>0$.
\begin{itemize}
  \item If $\alpha\leq1,$ then $u \in C^{0, \sigma}(\mathbb{R}^n)$ for any $\sigma<\alpha$.
  Moreover, $$\parallel u\parallel_{C^{0, \sigma}(\mathbb{R}^n)}\leq C
  (\parallel u\parallel_{L^\infty}+\parallel w\parallel_{L^\infty})$$
  for a constant $C$ depending only on $n, \alpha$ and $\sigma$.

  \item If $\alpha>1,$  then $u \in C^{1, \sigma}(\mathbb{R}^n)$ for any $\sigma
  <\alpha-1$.
  Moreover, $$\parallel u\parallel_{C^{1, \sigma}(\mathbb{R}^n)}\leq C
  (\parallel u\parallel_{L^\infty}+\parallel w\parallel_{L^\infty})$$
  for a constant $C$ depending only on $n, \alpha$ and $\sigma$.
\end{itemize}
\label{prop2.9}
\end{pro}

\begin{pro}\cite{Si}
Let $w=(-\lap)^{\alpha/2}u$. Assume that $w \in C^{0, \sigma}(\mathbb{R}^n)$ and
$u \in L^\infty$ for $\sigma \in (0,1]$ and $\alpha>0$.
\begin{itemize}
  \item If $\sigma+\alpha\leq1,$ then $u \in C^{0, \sigma+\alpha}(\mathbb{R}^n)$.
  Moreover, $$\parallel u\parallel_{C^{0, \sigma+\alpha}(\mathbb{R}^n)}\leq C
  (\parallel u\parallel_{L^\infty}+\parallel w\parallel_{C^{0, \sigma}})$$
  for a constant $C$ depending only on $n, \alpha$ and $\sigma$.

  \item If $\sigma+\alpha>1,$  then $u \in C^{1, \sigma+\alpha-1}(\mathbb{R}^n)$.
  Moreover, $$\parallel u\parallel_{C^{1, \sigma+\alpha-1}(\mathbb{R}^n)}\leq C
  (\parallel u\parallel_{L^\infty}+\parallel w\parallel_{C^{0, \sigma}})$$
  for a constant $C$ depending only on $n, \alpha$ and $\sigma$.
\end{itemize}
\label{prop2.8}
\end{pro}

Notice that $|v_k|\leq1$, applying the propositions above to (\ref{pe3}) gives:

 \begin{enumerate}
   \item if $\alpha\leq1$, then for all $0<\sigma<\alpha$,

 \emph{I.} $\parallel v_k \parallel_{C^{0,\alpha+\sigma}(\mathbb{R}^n)}<C_{n,\alpha,\sigma},$
        when $ \sigma+\alpha\leq1$.

   \emph{II.} $\parallel v_k \parallel_{C^{1,\alpha+\sigma-1}(\mathbb{R}^n)}<C_{n,\alpha,\sigma},
     $ when $\sigma+\alpha>1$.
 \item if $\alpha>1$, then for all $0<\sigma<\alpha-1$,

  \emph{III.} $\parallel v_k \parallel_{C^{l,\gamma}(\mathbb{R}^n)}<C_{n,\alpha,\sigma}, $
  where  $l$ is the integer part of $1+\alpha+\sigma$  and $\gamma=1+\alpha+\sigma-l$.
 \end{enumerate}

Part \emph{I} and part \emph{II} can be derived directly from
Proposition \ref{prop2.9} and Proposition \ref{prop2.8}.
To obtain part \emph{III}, we
first apply Proposition \ref{prop2.9} and it gives
$$\parallel v_k \parallel_{C^{1,\sigma}(\mathbb{R}^n)}
<C_{n,\alpha,\sigma}.$$
 Then we show that $ v_k$ has the regularity claimed above
 in a neighborhood of the origin. Because such reasoning can be done
 in a neighborhood of any point in $\mathbb{R}^n$,
 we conclude that $ v_k$ has the desired uniform regularity.

  Let $\varphi$ be a smooth
cutoff function such that $\varphi(x)\in [0,1]$
in $\mathbb{R}^n$, $supp\,\varphi \subset B_2$ and $\varphi(x)\equiv1$ in
$B_1$. Let $(-\lap)^{\alpha/2}v_k=g_k$. Define
$$v_k^o(x):=c_{n, -\alpha}\int_{\mathbb{R}^n}\frac{\varphi(y)g_k(y)}
{|x-y|^{n-\alpha}}dy=(-\lap)^{-\alpha/2}(\varphi g_k)(x).$$
Then $$(-\lap)^{\alpha/2}v_k^o(x)=g_k(x)=(-\lap)^{\alpha/2}v_k(x),\; x \in  B_1.$$
Hence $v_k^o-v_k$ is a smooth function in $B_{\frac{1}{2}}$. In addition,
its $C^{1,\sigma}-$norm can be controlled by the $L^\infty-$norms of
$v_k^o$ and $v_k$.

Given that
$$v_k^o=(-\lap)^{-\alpha/2}(\varphi g_k)=
(-\lap)^{1-\alpha/2}\circ(-\lap)^{-1}(\varphi g_k),$$
it then follows from the $C^{2,\sigma}-$estimates for the poisson equation
(see \cite{GT}) that $(-\lap)^{-1}(\varphi g_k)\in C^{3,\sigma}$. Also,
its norm is solely dependent of $\parallel v_k \parallel_{C^{1,\sigma}}$.
Thus applying Proposition \ref{prop2.7} gives
$v_k^o \in C^{l,\gamma}$ and $l>\alpha.$ This proves \emph{III}.

Due to the equi-continuity of $\{v_k\}$ in $\mathbb{R}^n$,
we are able to employ the \emph{Arzel$\grave{a}$-Ascoli theorem} and
 claim the existence of a converging sequence $\{v_{1m}\}$ in $B_1(0)$.
  For sure, one can find a subsequence of $\{v_{1m}\}$, denoted by $\{v_{2m}\}$, that converges in $B_2(0)$.
Then another subsequence of $\{v_{2m}\}$, denoted as $\{v_{3m}\}$, converges in $B_3(0)$. By induction, we get a chain of sequences
$$\{v_{jm}\}\supset\{v_{2m}\}\supset\{v_{3m}\}\cdot\cdot\cdot$$
such that $\{v_{jm}\}$ converges in $B_j(0)$.
Now form the diagonal sequence $\{v_{jj}\}$, whose $j$th term is the $j$th term in the
$j$th subsequence $\{v_{jm}\}$. Such $\{v_{jj}\}$ converges at all points in any $B_R(0)$. Thus we have
constructed a sequence of solutions that converges point-wise in $\mathbb{R}^n$.

To show the other part of (\ref{pe4}), we turn to the definition of the fractional Laplacian.
\begin{eqnarray*}
&&(-\lap)^{\alpha/2}v_k(x)\\
&=&\frac{1}{2}\int_{\mathbb{R}^n }\frac{2v_k(x)-v_k(x+y)-
v_k(x-y)}{|y|^{n+\alpha}}dy\\
&=& \frac{1}{2}\bigg[\int_{\mathbb{R}^n\backslash B_1(0) }\frac{2v_k(x)-v_k(x+y)-
v_k(x-y)}{|y|^{n+\alpha}}dy\\
&&+\int_{B_1(0)}\frac{2v_k(x)-v_k(x+y)-
v_k(x-y)}{|y|^{n+\alpha}}dy\Bigg]\\
&=& \frac{1}{2}(I_1+I_2).
\end{eqnarray*}
Henceforth we use $C$ to denote positive constants whose values may vary from line to line.
Doing Taylor expansion near $x$ to $v_k$ in $I_2$, the equi-continuity of $v_k$ gives
\begin{eqnarray*}
|I_2|&\leq&\int_{B_1(0) }\frac{|2v_k(x)-v_k(x+y)-
v_k(x-y)|}{|y|^{n+\alpha}}dy\\
&\leq& \int_{B_1(0) }\frac{C|y|^{\alpha+\epsilon}}{|y|^{n+\alpha}}dy\\
&=&\int_{B_1(0) }\frac{C}{|y|^{n-\epsilon}}dy\\
&<&\infty.
\end{eqnarray*}
Meanwhile, since $|v_k|\leq1$,
\begin{eqnarray*}
|I_1|&\leq&\int_{\mathbb{R}^n\backslash B_1(0)}\frac{|2v_k(x)-v_k(x+y)-
v_k(x-y)|}{|y|^{n+\alpha}}dy\\
&\leq&\int_{\mathbb{R}^n\backslash B_1(0)}\frac{4C}{1+|y|^{n+\alpha}}dy\\
&<&\infty.
\end{eqnarray*}
Therefore,
\begin{eqnarray*}
&&\int_{\mathbb{R}^n }\left|\frac{2v_k(x)-v_k(x+y)-
v_k(x-y)}{|y|^{n+\alpha}}\right|dy\\
&\leq&\int_{\mathbb{R}^n}\frac{C}{1+|y|^{n+\alpha}}\left(1+
\frac{1}{|y|^{n-\epsilon}}\right)dy\\
&<& \infty.
\end{eqnarray*}
It then follows from the \emph{Lebesgue's dominated convergence theorem} that
\begin{eqnarray*}
&&\lim_{k \ra \infty}(-\lap)^{\alpha/2}v_k(x)\\
&=&\lim_{k \ra \infty}\frac{1}{2}\int_{\mathbb{R}^n }\frac{2v_k(x)-v_k(x+y)-
v_k(x-y)}{|y|^{n+\alpha}}dy\\
&=&\frac{1}{2}\int_{\mathbb{R}^n }\lim_{k \ra \infty}\frac{2v_k(x)-v_k(x+y)-
v_k(x-y)}{|y|^{n+\alpha}}dy\\
&=&(-\lap)^{\alpha/2}v(x).
\end{eqnarray*}
This proves (\ref{pe4}) and (\ref{pe5}).
\vspace{12pt}

\emph{Case ii.} $\underset{k \ra \infty}{\lim}\frac{d_k}{\lambda_k} = \;C\;>\;0$.
\vspace{12pt}

In this case,
 $$ \Omega_k \ra R^n_{+C}:=\{x_n\geq -C\mid x \in R^n \} \mbox{ as }k \ra \infty.$$
 Similar to \emph{Case i}, here we're able to establish the existence of a function $v$ such that as $k \ra \infty$,
\be
v_k(x) \ra v(x) \mbox{ and }
 (-\lap)^{\alpha/2}v_k(x)\ra(-\lap)^{\alpha/2}v(x),\label{pe6}
\ee
and thus
\be
(-\lap)^{\alpha/2}v(x)=v^p(x), \quad x \in  \mathbb{R}^n_{+C}. \label{pe7}
\ee

It's known that (\ref{pe7})
has no positive solution (see ). Meanwhile, it follows from (\ref{pe2.1}) that
$$v(0)=\lim_{k \ra \infty} v_k(0)=1.$$
This is a contradiction. Hence (\ref{pe2}) is true. Next we prove
(\ref{pe6}) and (\ref{pe7}).

Let $D_1=B_1(0)\cap \{x_n\geq0\}$. Through an argument similar to that in \emph{Case i}, we have the equicontinuity of $v_k$ and
$$|v_k(x)-v_k(y)|\leq C_{D_1}|x-y|^{\alpha+\epsilon},
\quad x,y\in D_1, \quad0<\epsilon<1.$$
Together with $|v_k(x)|\leq1$,
we are able to show the existence of a converging subsequence $\{v_{1k}\}$ of $\{v_k\}$
such that for each given $x\in D_1$,
$$v_{1k}(x) \ra v(x)\mbox{ and }(-\lap)^{\alpha/2}v_{1k}(x) \ra
(-\lap)^{\alpha/2}v(x).$$
Let $D_2=B_2(0)\cap \{x_n\geq-\frac{C}{2}\}$. We can still have
$$|v_k(x)-v_k(y)|\leq C_{D_2}|x-y|^{\alpha+\epsilon},\quad x,y\in D_2.$$
This again gives a converging
subsequence $\{v_{2k}\}$ of $\{v_{1k}\}$ such that for each given $x\in D_2$,
$$v_{2k}(x) \ra v(x)\mbox{ and }
(-\lap)^{\alpha/2}v_{2k}(x) \ra (-\lap)^{\alpha/2}v(x).$$
Let $D_3=B_3(0)\cap \{x_n\geq-\frac{2C}{3}\}$. There exists a
subsequence $\{v_{3k}\}$ of $\{v_{2k}\}$ such that for each given $x\in D_3$,
$$v_{3k}(x) \ra v(x)\mbox{ and }
(-\lap)^{\alpha/2}v_{3k}(x) \ra (-\lap)^{\alpha/2}v(x).$$ Repeating the process
above for $m$ times and it gives a pointwise converging subsequence
$\{v_{mk}\}$ in $D_m=B_m(0)\cap \{x_n\geq-\frac{(m-1)}{m}C\}$. Now we take the diagonal sequence $\{v_{mm}\}$ with $v_{mm}(x)$ being the $m$th term in $\{v_{mk}\}$. It's easy to see that for a fixed $x_0 \in \mathbb{R}^n_{+C}$, there exists
some $m_{x_0}$ such that for $m>m_{x_0}+1$,
$$v_{mm}(x) \ra v(x)\mbox{ and }
(-\lap)^{\alpha/2}v_{mm}(x) \ra (-\lap)^{\alpha/2}v(x).$$
Thus verifies (\ref{pe6}) and (\ref{pe7}).
\vspace{12pt}

\emph{Case iii.} $\underset{k \ra \infty}{\lim}\frac{d_k}{\lambda_k} = \;0$.
\vspace{12pt}

In this case, there exists a point $x^o \in \partial\Omega$ such that
$$x^k \:\ra\: x^o, \quad k \:\ra \:\infty.$$
Let $\tilde{x^o}=\frac{x^o-x^k}{\la_k}$. Obviously,
$$|\tilde{x^o}| \:\ra \:0, \quad k\: \ra\: \infty.$$
Then we will show that $v_k$ is uniformly H$\ddot{o}$lder continuous near $\tilde{x^o}$, i.e.
\be
|v_k(x)-v_k(\tilde{x^o})|\leq C|x-\tilde{x^o}|^{\alpha/2}.\label{pe9}
\ee
We postpone the proof of (\ref{pe9}) for a moment.
Notice that
\begin{equation}\label{pe10}
v_k(0)-v_k(\tilde{x^o})=1.
\end{equation}
On the other hand, it follows from (\ref{pe9}) that
\begin{equation}\label{pe11}
  |v_k(x)-v_k(\tilde{x^o})| \ra\, 0, \mbox{ as } x \ra\, 0.
\end{equation}
This contradicts (\ref{pe10}). Therefore \emph{Case iii} will not happen.

To prove (\ref{pe9}), we introduce an auxiliary
function $\varphi$ such that for all $x$ near $\tilde{x^o}$, it holds that
\be\label{pe13}
|v_k(x)-v_k(\tilde{x^o})|\leq
|\varphi(x)-\varphi(\tilde{x^o})|\leq C|x-\tilde{x^o}|^{\alpha/2}.
\ee

Let \begin{equation}
\psi_1(x)=\left\{\begin{array}{ll}
C_{n, \alpha}(1-|x-z|^2)^{\alpha/2},  &x \in B_1(z),\\
0,&x \in \overline{B_1^c(z)},
\end{array}
\right.
\end{equation}
where $z=\frac{\tilde{x^o}}{|\tilde{x^o}|}(1+\tilde{x^o})$  and
$B_1^c(z)=\mathbb{R}^n\backslash B_1(z)$.
To simplify the calculation, without any loss of generality,
we set $z$ to be the origin, thus
\begin{equation}
\psi_1(x)=\left\{\begin{array}{ll}
C_{n, \alpha}(1-|x|^2)^{\alpha/2},  &x \in B_1(0),\\
0,&x \in \overline{B_1^c(0)}.
\end{array}
\right.
\end{equation}
Let
\begin{equation*}
  \psi_2(x)=\frac{1}{|x|^{n-\alpha}}\psi_1(\frac{x}{|x|^2})
\end{equation*}
 be the Kelvin transform of $\psi_1(x)$. Then
\begin{equation*}
 (-\lap)^{\alpha/2}\psi_2(x)=\frac{C_{n, \alpha}}{|x|^{n+\alpha}}, \quad x \in B_1^c(0).
\end{equation*}
Let $\xi(x)$ be a smooth cutoff function such that $\xi(x) \in [0, 1] $ in
$R^n$ and $\xi(x)=0$ in $\overline{B_1(0)}$ and $\xi(x)=1$ in $B_3^c(0)$.
Then we have
\begin{eqnarray*}
(-\lap)^{\alpha/2}\xi(x)
  &=& \frac{1}{2}\int_{R^n }
  \frac{2\xi(x)-\xi(x+y)-\xi(x-y)}{|y|^{n+\alpha}}dy \\
   &\geq&-C.
\end{eqnarray*}

For $k>0$, let $$\varphi(x)=k\psi_2(x)+\xi(x),$$
 and $D=(B_3(0)\backslash B_1(0)) \cap \Omega_k.$
 For $k$ sufficiently large and a fixed $x \in D$,
 \begin{eqnarray*}
  (-\lap)^{\alpha/2}\varphi(x) &=&
  k(-\lap)^{\alpha/2}\psi_2(x) +(-\lap)^{\alpha/2}\xi(x) \\
    &\geq&\frac{kC_{n, \alpha}}{|x|^{n+\alpha}}-C \\
    &\geq&1.
 \end{eqnarray*}
 Thus
 \begin{equation}\label{pe12}
\left\{\begin{array}{ll}
 (-\lap)^{\alpha/2}(\varphi-v_k)\geq0,  &x \in D,\\
\varphi(x)-v_k(x)>0, &x \in D^c.
\end{array}
\right.
\end{equation}
Applying the \emph{maximum principle} (see \cite{Si}) to (\ref{pe12}), it gives
\begin{equation*}
  \varphi(x)\geq v_k(x), \quad x \in D.
\end{equation*}

To show that $\varphi(x)$ is H$\ddot{o}$lder continuous near $\tilde{x^o}$ in $D$, it suffices
to show that
$\psi_2(x)$ has such H$\ddot{o}$lder continuity. Indeed,
\begin{eqnarray*}
 \psi_2(x)-\psi_2(\tilde{x^o})
 &=&\frac{1}{|x|^n}(|x|-1)^{\alpha/2}(|x|+1)^{\alpha/2}\\
  &\leq& C\, (|x|-1)^{\alpha/2}.
\end{eqnarray*}
Hence
\be
0\leq v_k(x)-v_k(\tilde{x^o})= v_k(x)\leq \varphi(x)=
\varphi(x)-\varphi(\tilde{x^o})\leq C|x-\tilde{x^o}|^{\alpha/2}.
\ee
This proves (\ref{pe13}) and completes the proof of \emph{Case iii}.

\bigskip

{\em Authors' Addresses and E-mails:}
\medskip

Wenxiong Chen

Department of Mathematical Sciences

Yeshiva University

New York, NY, 10033 USA

wchen@yu.edu
\medskip

Congming Li

Department of Mathematics, INS and MOE-LSC

Shanghai Jiao Tong University

Shanghai, 200240, China, and

Department of Applied Mathematics

University of Colorado,

Boulder CO USA

cli@clorado.edu
\medskip

Yan Li

Department of Mathematical Sciences

Yeshiva University

New York, NY, 10033 USA

yali3@mail.yu.edu

\end{document}